\documentclass{article}

\usepackage[dcucite]{harvard}

\usepackage{amssymb,amsmath}
\usepackage[latin1]{inputenc}

\newtheorem{theorem}{Theorem}

\newtheorem{definition}[theorem]{Definition}

\newtheorem{lemma}[theorem]{Lemma}
\newtheorem{remark}[theorem]{Remark}
\newenvironment{proof}[1][Proof]{\textbf{#1.} }{\ \rule{0.5em}{0.5em} \bigskip}

\newenvironment{keywords}{\begin{center}
\begin{minipage}[c]{11cm} {\bf Keywords:}} {\end{minipage}
\end{center}}
\newenvironment{msc}{\begin{center}
\begin{minipage}[c]{11cm} {\bf 2000 Mathematics Subject Classification:}} {\end{minipage}
\end{center} \bigskip}

\begin{document}

\title{Existence and uniqueness of solutions
for a nonlocal parabolic thermistor-type problem\footnote{To be
presented at the \emph{13th IFAC Workshop on Control Applications of
Optimisation}, 26-28 April 2006, Paris -- Cachan, France. Accepted
(19-12-2005) for the Proceedings, IFAC publication, Elsevier Ltd,
Oxford, UK. Research Report CM05/I-55.}}

\author{Abderrahmane El Hachimi${}^{\dag}$\\
\texttt{elhachimi@ucd.ac.ma} \and
Moulay Rchid Sidi Ammi${}^{\ddag}$\\
\texttt{sidiammi@mat.ua.pt} \and
Delfim F. M. Torres${}^{\ddag}$\\
\texttt{delfim@mat.ua.pt}}

\date{${}^{\dag}$UFR Math\'ematiques Appliqu\'ees et Industrielles\\
Facult\'{e} des Sciences, Universit\'e Chouaib Doukkali\\
B.P. 20, El Jadida, Maroc\\
[0.3cm]
${}^{\ddag}$Department of Mathematics\\
University of Aveiro\\
3810-193 Aveiro, Portugal}

\maketitle


\begin{abstract}
In this paper we prove existence and uniqueness of solutions to the
nonlocal parabolic problem
$$
  \frac{\partial u}{\partial t}-\triangle_{p} u
  = \lambda \frac{f(u)}{\left( \int_{\Omega} f(u)\, dx \right)^2} , \quad
   \mbox{in } \Omega \times ]0,T[ \, ,
$$
which generalizes the electric heating problem of a conducting body.
\end{abstract}


\begin{msc}
35K55, 35B41, 80A20.
\end{msc}


\section{Introduction}

In this paper we study the existence and uniqueness of
bounded solutions for the following nonlocal parabolic problem:
\begin{gather}
\label{11}
  \frac{\partial u}{\partial t}-\triangle_{p} u = \lambda
  \frac{f(u)}{\left( \int_{\Omega} f(u)\, dx \right)^2} \, , \quad
   \mbox{in } \Omega \times ]0,T[ \, , \notag \\
   u = 0 \quad \mbox{on }  \partial \Omega \times ]0,T[\, ,  \\
   u(0)= u_0  \quad \mbox{in }  \Omega \, , \notag
\end{gather}
where $\Delta_{p}= div(-|\nabla u|^{p-2}\nabla u)$; $p \ge 2$;
$T>0$; $\Omega \subset \mathbb{R}^{N}$, $N\geq 1$, is a regular
bounded domain; $\lambda$ a positive parameter; and $f$ a function
from $\mathbb{R}$ to $\mathbb{R}$ with prescribed conditions.

For $p=2$, $\Delta_{p}$ is reduced to the usual Laplacian operator, and
problem (\ref{11}) serves as a model for the well-known and
important thermistor problem, where $u$ is the temperature inside
a conductor -- see \textrm{e.g.} \cite{lac1,lac2,bl,tza}. This
problem is very important in industry and engineering applications,
and has attracted attention in the literature over the last decade,
from both the experimental and theoretical point of views: see
\cite{ac,al,es1,MR1903000,es2,MR2108729} and references therein.

Our main result is a proof of the global existence and uniqueness of
solutions of problem (\ref{11}). The result is a generalization of
\cite{lac1,lac2,tza,es2} to the general $p$-Laplacian case, $p \ge
2$. For the particular case $p=2$, the result is obtained in
\cite{es2}, but under somehow less restrictive assumptions on the
data of the problem: Theorem~\ref{thm:mainResult:SA} does not impose
restrictions on $\alpha$, while in \cite{es2} it is assumed that
$\alpha < \frac{4}{N-2}$, $N>2$.


\section{Existence and uniqueness}

The definition of solution for problem (\ref{11})
is understood in the standard way.

\begin{definition}
We say that $u$ is a solution of (\ref{11}) if, and only if,
$$
 u \in L^{\infty}(\tau ,+\infty ,W_{0}^{1,p}(\Omega)\cap
L^{\infty}(\Omega))
$$
with $\frac{\partial u}{\partial t}\in L^{2}(\tau ,+\infty
,L^{p'}(\Omega) )$ for any  $\tau >0$, and the following equation is
satisfied for all $ \phi \in C^{\infty}((0,\infty),\Omega )$:
\begin{equation*}
\int_{0}^{T}\int_{\Omega} u \frac{\partial }{\partial t}\phi
-|\nabla u|^{p-2}
 \nabla u \nabla \phi \, dxdt
 = \int_{0}^{T}\left(  \frac{\lambda}{ \left(\int_{\Omega} f(u)\, dx
\right)^{2}}
 \int_{\Omega} f(u) \phi dx\right)dt \, .
\end{equation*}
\end{definition}

The main result of the paper is as follows.

\begin{theorem}
\label{thm:mainResult:SA} Let the hypotheses (H1) and (H2) be
satisfied:

 \begin{itemize}

\item[(H1)] $f: \mathbb{R} \rightarrow \mathbb{R} $ is a locally
Lipschitzian function;

\item[(H2)] There exist positive constants $c_{1}$, $c_{2}$ and
$\alpha$ such that for all $\xi \in \mathbb{R}$
  $$
  \sigma \leq f(\xi )\leq c_{1}| \xi|^{\alpha +1}+c_{2} \, .
  $$

\end{itemize}

Further, assume that $u_{0}\in L^{k_{0}+2}(\Omega)$ with
  \begin{equation} \label{21}
     k_{0} \geq \max \left( 0,\frac{N (\alpha +2-p)}{p}-2 \right).
   \end{equation}
Then, there exists a constant $d_{0}>0$ such that if
   $\|u_{0}\|_{k_{0}+2}<d_{0}$,
   the problem (\ref{11}) admits a solution $u$ verifying
\begin{gather*}
   u \in L^{\infty}(\tau ,+\infty ,L^{k_{0}+2}(\Omega)) \, , \\
  |u|^{\gamma}u \in L^{\infty}(\tau ,+\infty ,W_{0}^{1,p}(\Omega)),
    \mbox{ with } \gamma = \frac{k_{0}}{p} \, ,
\end{gather*}
for all $\tau >0$. Moreover, if $u_{0}\in L^{\infty}(\Omega ),$ then
$u \in L^{\infty}(\tau ,+\infty ,L^{\infty }(\Omega))$ and $u$ is
unique.
\end{theorem}

\begin{remark}
A value for $d_{0}$ is given explicitly in the proof of
Theorem~\ref{thm:mainResult:SA} -- \textrm{cf.} \eqref{28}.
\end{remark}


\section{Proof of Theorem~\ref{thm:mainResult:SA}}

The existence is proved by the Faedo-Galerkin method.


\subsection{Existence}

Let $w_{1}, \ldots, w_{m}, \ldots$ be a complete sequence of
linearly independent elements of $H_{0}^{1}(\Omega)$. For each $m$,
we define an approximate solution
\begin{equation*}
u_{m}(t)= \sum_{j=1}^{m} g_{jm}(t)w_{j}
\end{equation*}
of (\ref{11}), where $g_{jm}$ are solutions of the following
system of ordinary differential equations:
\begin{equation}
\label{22}
\langle u_{m}',w_{j}\rangle +(u_{m},w_{j})
= \frac{\lambda }
 { \left(\int_{\Omega} f(u_{m})\, dx \right)^{2} } \,
 \langle f(u_{m}),w_{j} \rangle \, ,
\end{equation}
$j = 1,\ldots,m$, with the initial condition
\begin{equation}
\label{eq:ic} u_{m}(0)=u_{om} \, ,
\end{equation}
$u_{om}$ being the orthogonal projection in
 $H_{0}^{1}(\Omega)$ of $u_{0}$ on the space spanned by
 $w_{1}, \ldots, w_{m}$. The initial-value problem
 (\ref{22})-(\ref{eq:ic}) is equivalent to a
linear $m$-dimensional ordinary differential equation for the
$g_{jm}$. The existence and uniqueness of the $g_{jm}$ on a maximal
interval $[0,t_{m}[$ is obvious. We obtain the existence of a solution
$u$ for our problem \eqref{11} passing to the limit, as $m
\rightarrow \infty$. For that we need to derive a priori estimates
on $u_{m}$ which guarantee that $t_{m}=T$. This is done by
Lemma~\ref{lem:estimates}. In order to prove it, we employ
an inequality due to Ghidaglia.

\begin{lemma}[Ghidaglia inequality]
\label{eq:Ghidaglia} Let $y$ be a positive absolutely continuous
function on $(0,+\infty)$ which satisfies
$$ y' + \gamma y^{\nu} \leq \delta \, , $$
with $\nu>1, \gamma >0$ and $\delta \geq 0$. Then,
$$ y(t) \leq \left(\frac{\delta}{\gamma}\right)^{\frac{1}{\nu}}
+ \left(\gamma (\nu -1)t\right)^{-\frac{1}{(\nu -1)}} \, ,$$ for all
$t\geq 0$,
\end{lemma}

Proof of Lemma~\ref{eq:Ghidaglia} can be found in \cite{tem}.

 \begin{lemma}
 \label{lem:estimates}
  For any $ \tau >0$, there exist constants $c_{3}(\tau)$ and
$c_{4}(\tau)$
  such that for all $t \geq \tau$
  \begin{equation}\label{23}
  \| u_{m}(t)\|_{k_{0}+2} \leq  c_{3}(\tau ) \, ,
  \end{equation}
  \begin{equation}\label{24}
  \| u_{m}(t)\|_{\infty } \leq  c_{4}(\tau ) \,  .
 \end{equation}
 \end{lemma}

\begin{remark}
 Throughout the paper we denote by $c_{i}$ different positive
constants,
 which depend on the data of the problem, but not on $m$.
\end{remark}

\begin{proof}
Multiplying the equation (\ref{22}) by $|u_{m}|^{k} g_{jm}$,
integrating on $\Omega $, summing up for $j=1,\ldots,m$ and using
(H1)-(H2), yields
 \begin{equation}\label{25}
 \frac{1}{k+2}\frac{d}{dt}\| u_{m}\|_{k+2}^{k+2}
+\frac{p^{p}}{(k+p)^{p}}
     \| \nabla (|u_{m}|^{\frac{k}{p}} u_{m})\|_{p}^{p}
     \leq c_{5} \| u_{m}\|_{k+\alpha +2}^{k+\alpha +2}+ c_{6} \, .
\end{equation}
 By using condition \eqref{21} on $k_{0}$ and well-known Sobolev's
and Gagliardo-Nirenberg's inequalities, we obtain
\begin{equation*}
     \left(c_{7} \| u_{m}\|_{k_{0}+2}^{\alpha }-
     \frac{4}{(k_{0}+p)^{p}}\right)\| \nabla |u_{m}|^{\gamma
}u_{m}\|_{p}^{p}
     + c_{6} \\
\ge \frac{1}{k_{0}+2}\frac{d}{dt}\| u_{m}\|_{k_{0}+2}^{k_{0}+2} \, .
\end{equation*}

Using the compatibility condition on $u_{0}$
\begin{equation} \label{28}
        \| u_{0}\|_{k_{0}+2} <
         \left( \frac{4}{c_{7}(k_{0}+p)^{p}}\right)^{\frac{1}{\alpha}}
         = d_{0} \, ,
\end{equation}
and the continuity of $u_{m}$, there exists a small $\tau >0$ such
that
 \begin{equation}
 \label{210}
  \frac{1}{k_{0}+2}\frac{d}{dt}\| u_{m}\|_{k_{0}+2}^{k_{0}+2}
         + c_{8}\| \nabla (|u_{m}|^{\gamma }u_{m})\|_{p}^{p}
          \leq c_{6}
    \end{equation}
for all $0<t<\tau$. Setting
$$y_{k_{0}}(t)=\|u_{m}\|_{k_{0}+2}^{k_{0}+2}$$ and using the
Poincar\'e and Holder inequalities on the left side of (\ref{210}),
    there exist two
    constants $\gamma >0$ and $\delta >0$ such that
    $$
    \frac{dy_{k_{0}}}{dt}+ \gamma
    y_{k_{0}}^{\frac{k_{0}+p}{k_{0}+2}}\leq \delta
$$
for all $0<t<\tau$. Note that for $p>2$ we have
$\frac{k_{0}+p}{k_{0}+2} >1$. Estimate \eqref{23} follows from
Lemma~\ref{eq:Ghidaglia}.

The proof of \eqref{24} is similar to the proof of inequality (2.4)
in \cite{es2}, and is given here for completeness.
By using Holder's inequality, we get
   \begin{equation} \label{211}
           \| u_{m}\|_{k+\alpha +2}^{k+\alpha +2}\leq c_{9}
            \| u_{m}\|_{k+2}^{\theta _{1}}   \| u_{m}\|_{k_{0}+2}^{\theta _{2}}
             \| u_{m}\|_{q}^{\theta _{3}} \, ,
      \end{equation}
 with $\theta _{1}, \theta _{2}$ and $\theta _{3}$ satisfying
 $$
  \frac{\theta _{1}}{k+2}+ \frac{\theta _{2}}{k_{0}+2}+
  \frac{\theta _{3}}{q} =1 $$
   and $$
   \theta _{1}+\theta _{2}+\theta _{3}=k+\alpha +2.
 $$
  Moreover, we require
   $$
  \frac{\theta _{1}}{k+2}+ \frac{\theta _{3}}{p(\gamma +1)}=1 .
  $$
   Using the boundedness of $\| u_{m}\|_{k_{0}+2}$, the choice of $q$,
   Sobolev and Young's inequalities and relation  \eqref{211}, we derive
\begin{multline*}
c_{5} \| u_{m}\|_{k+\alpha +2}^{k+\alpha +2}\leq c_{10}
\| u_{m}\|_{k+2}^{\theta _{1}} \|
\nabla |u_{m}|^{\gamma}u_{m}\|_{p}^{\frac{\theta_{3}}{\gamma +1}} \\
\leq c_{11}(k+2)^{\theta _{4}} \| u_{m}\|_{k+2}^{k+2}
+ \frac{p^{p}}{2(k+p)^{p}}
\| \nabla |u_{m}|^{\gamma }u_{m}\|_{2}^{2} \, ,
\end{multline*}
      where $\theta _{4}$ is some positive constant.
      Hence  $(\ref{25})$ becomes
\begin{equation*}
        \frac{1}{k+2}\frac{d}{dt}\|u_{m}\|_{k+2}^{k+2} +\frac{c_{12}}{(k+p)^{p}}
        \| \nabla |u_{m}|^{\gamma }u_{m}\|_{p}^{p}
        \leq
         c_{13}(k+p)^{\theta _{4}} \| u_{m}\|_{k+2}^{k+2} +c_{6}.
\end{equation*}
Therefore, by applying Lemma~4 of \cite{jf}, we conclude \eqref{24}.
\end{proof}

Multiplying the jth equation of system (\ref{22}) by $g_{jm}(t)$,
summing these equations for $j=1,\ldots, m$ and integrating with
respect to the time variable, we deduce the existence of a
subsequence of $u_{m}$ such that
\begin{gather*}
u_{m} \to u  \mbox{ weak star in } L^{\infty}(0,T;L^{2}(\Omega)) \,
,
\\
u_{m} \to u  \mbox{ weak in } L^{2 }(0,T;W^{1,p}_{0}(\Omega)) \, , \\
u_{mt} \to u_{t}  \mbox{ weak in } L^{2}(0,T;W^{-1,p'}(\Omega)) \, , \\
u_{m} \to u  \mbox{ strongly in } L^{p }(0,T;L^{p}(\Omega)) \, .
\end{gather*}
Standard compactness and monotonicity arguments allow us to assert
that $u$ is a solution of problem (\ref{11}).


\subsection{Uniqueness}

Let $u_{1}$ and $u_{2}$ be two weak solutions of problem
(\ref{11}), and define $w=u_{1}-u_{2}$. Subtracting the equations
verified by $u_{1}$ and $u_{2}$, we obtain:
\begin{multline*}
\frac{dw}{dt}-(\triangle_{p} u_{2} -\triangle_{p} u_{1})
= \frac{\lambda \left( f(u_{1})-f(u_{2}) \right)}{\left(
\int_{\Omega}f(u_{1})\, dx \right)^{2}} \\
+\lambda \frac{\left( \int_{\Omega }
f(u_{2})-f(u_{1})\,dx \right)
   \left( \int_{\Omega }f(u_{2})+f(u_{1})\,dx \right)}
  {\left( \int_{\Omega} f(u_{1})\, dx \right)^{2}
  \left( \int_{\Omega} f(u_{2})\, dx \right)^{2}} f(u_{2}) \, .
\end{multline*}
Taking the inner product of  last equation by $w$ and using (H1),
(\ref{24}), and the monotonicity of the $p$-Laplacian, we get
  $$
  \frac{1}{2}\frac{d}{dt}\| w(t)\|^{2}_{2}\leq c_{14} \|w(t)\|^{2}_{2}
\, ,
  $$
which implies that $w=0$. Hence, the solution is unique.


\section{Absorbing sets and attractors}

We denote by $\{ S(t), t\geq 0\}$ the continuous semi-group
generated by
  $(\ref{11})$ and defined by
   $$
  \begin{array}{rcl}
  S(t): L^{\infty }(\Omega ) &\rightarrow & L^{\infty }(\Omega ) \\
  u_{0} & \rightarrow & S(t)u_{0}=u(t,.).
  \end{array}
  $$

Using the techniques of R.~Temam \cite{tem}, we prove
existence of attractors.

\begin{theorem}\label{thm31}
 The semigroup $S(t)$, associated with the problem (\ref{11}),
 possesses a maximal attractor $A$ which is bounded in
 $W_{0}^{1,p}(\Omega)$, compact and connected  in
 $L^{\infty}(\Omega)$.
 \end{theorem}

\begin{proof}
Inequality (\ref{24}) implies that there exists an absorbing set in
$L^{k}(\Omega), 1\leq k \leq \infty$. We now prove the existence of
an absorbing set in $W_{0}^{1,p}(\Omega)$ and the uniform
compactness of the semigroup $S(t)$. For this purpose, multiplying
(\ref{22}) by $g'_{jm}(t)$, summing up from $j=1$ to $m$, integrating
over $\Omega$ and using Holder inequality, one obtains that
\begin{equation*}
 \left\|\frac{\partial u_{m}}{\partial t}\right\|_{2}^{2}
 +\frac{1}{p}\frac{\partial}{\partial
t}\left\|u_{m}\right\|_{W_{0}^{1,p}(\Omega)}^{p} \\
\leq c_{15}\int f(u_{m})\frac{\partial u_{m}}{\partial t} \\
\leq  c_{16}(\tau)+ \frac{1}{2}\left\|\frac{\partial u_{m}}{\partial
t}\right\|_{2}^{2} \, .
\end{equation*}
We deduce that for all $t \geq \tau$
\begin{equation*}
\left\|\frac{\partial u_{m}}{\partial t}\right\|_{2}^{2}
+\frac{\partial}{\partial
t}\left\|u_{m}\right\|_{W_{0}^{1,p}(\Omega)}^{p} \leq c_{17}(\tau)
\, .
\end{equation*}
 Hence,
 \begin{equation} \label{214}
\frac{\partial}{\partial
t}\left\|u_{m}\right\|_{W_{0}^{1,p}(\Omega)}^{p} \leq c_{17}(\tau)
\, , \quad \forall t \geq \tau \, .
 \end{equation}
Multiplying (\ref{22}) by $g_{jm}(t)$ we also have
\begin{equation} \label{215}
\begin{array}{ll}
\frac{1}{2}\frac{\partial}{\partial
t}\|u_{m}\|_{2}^{2}+\|u_{m}\|_{W_{0}^{1,p}(\Omega)}^{p} &\leq
c_{18}\int |f(u_{m}) u_{m}| \\ & \leq c_{19}(\tau).
\end{array}
 \end{equation}
After integrating in $t$, we infer from the last
equation (\ref{215}) that
 \begin{equation} \label{216}
\int_{t}^{t+\tau}\|u_{m}\|_{W_{0}^{1,p}(\Omega)}^{p} \leq
c_{20}(\tau) \quad \forall t \geq \tau .
 \end{equation}
Using (\ref{214})-(\ref{216}), we can apply  the uniform Gronwall Lemma
\cite[p.~89]{tem}, and by the lower semi-continuity of the norm, we
conclude that
$$
\|u_{m}\|_{W_{0}^{1,p}(\Omega)}^{p} \leq c_{21}(\tau) \, , \quad
\forall t \geq \tau \, .
$$
It follows that the ball $B(0,c_{21}(\tau))$ of
$W_{0}^{1,p}(\Omega)$, centered at $0$ and with radius
$c_{21}(\tau)$, is absorbing in $W_{0}^{1,p}(\Omega)$. The
assumption of Theorem I.1.1 in \cite[p.~23]{tem} is satisfied, and
the proof of Theorem \ref{thm31} is complete.
\end{proof}


\section{Conclusions and future work}

In this paper we prove existence and uniqueness for a
$p$-Laplacian nonlinear system of partial differential equations of
parabolic type, $p \ge 2$. For $p = 2$ the problem is a model of the heat
diffusion produced by the Joule effect in an electric conductor, and
we recover the previously known existence, boundedness, and
uniqueness results found in the literature for the thermistor
problem.

In a forthcoming work we will investigate the possibility to prove
more regularity results of the solution of the problem, by imposing
more restrictive assumptions on the data. The question is nontrivial
due to the nonlinear nature of the problem, as shown in
\cite{MR2079805} for $p = 2$.


\section*{Acknowledgment}

M. R. Sidi Ammi acknowledges the support of FCT (\emph{The
Portuguese Foundation for Science and Technology}), fellowship
SFRH/BPD/20934/2004.



\normalsize

\end{document}